\documentstyle[12pt]{article}
\def\Bbb R{{\rm \bf R}}
\def\proclaim#1{\vskip2mm{\bf #1}\em}
\def\endproclaim{\em \vskip2mm}
\def\tag#1{\eqno(#1)}
\def\gathered{\begin{array}{c}}
\def\endgathered{\end{array}}
\def\text{\mbox}

\begin{document}

\title {Prescribing a heat flux coming from a wave equation}
\author{Masaru IKEHATA\footnote{
Laboratory of Mathematics,
Graduate School of Engineering,
Hiroshima University,
Higashihiroshima 739-8527, JAPAN}}
%\date{}
\maketitle

\begin{abstract}
What happens when one prescribes a heat flux which is proportional to the Neumann data of
a solution of the wave equation in the whole space on the surface of a heat conductive body?  
It is shown that there is a difference in the asymptotic behaviour of the indicator function
in the most recent version of {\it the time domain enclosure method}, which aims at extracting
information about an unknown cavity embedded in the body.

\noindent
AMS: 35R30, 35K05, 35L05

\end{abstract}

%{\bf SendingWavelikeHeatFlux.tex}

\section{Introduction}

Let $\Omega$ be a bounded domain of $\Bbb R^3$ with $C^2$-boundary.
Le $D$ be a nonempty bounded open subset of $\Omega$ with $C^2$-boundary such that $\overline D\subset\Omega$ and 
$\Omega\setminus\overline D$ is connected.
Let $0<T<\infty$.

Given $f=f(x,t), (x,t)\in\partial\Omega\times\,]0,\,T[$, which belongs to $L^2(0,T;H^{-1/2}(\partial\Omega))$,
let $u=u_f(x,t)$, with $(x,t)\in\,(\Omega\setminus\overline D)\times\,]0,\,T[$, 
denote the weak solution of the following initial boundary value problem for the heat equation:
$$\displaystyle
\left\{
\begin{array}{ll}
\displaystyle
(\partial_t-\Delta) u=0 & \text{in}\,(\Omega\setminus\overline D)\times\,]0,\,T[,\\
\\
\displaystyle
u(x,0)=0 & \text{in}\,\Omega\setminus\overline D,
\\
\\
\displaystyle
\frac{\partial u}{\partial\nu}=0 & \text{on}\,\partial D\times\,]0,\,T[,\\
\\
\displaystyle
\frac{\partial u}{\partial\nu}=f(x,t) & \text{on}\,\partial\Omega\times\,]0,\,T[.
\end{array}
\right.
\tag {1.1}
$$
We use the same symbol $\nu$ to denote both the outer unit normal vectors of both $\partial D$
and $\partial\Omega$. 
The solution class is the same one as in \cite{IK1} which employs the weak solution in \cite{DL}.
The $u$ satisfies, for a positive constant $C_T$ independent of $f$,
$$\displaystyle
\Vert u\Vert_{L^2(0,T;H^1(\Omega\setminus\overline D))}
+\Vert\partial_tu\Vert_{L^2(0,T;H^1(\Omega\setminus\overline D)')}
+\Vert u(\,\cdot\,,T)\Vert_{L^2(\Omega\setminus\overline D)}
\le C_T\Vert f\Vert_{L^2(0,T;H^{-1/2}(\partial\Omega))}.
$$
See Section 2.1 in \cite{IK1} for more information about the direct problem.

This paper is concerned with the inverse obstacle problem described below.

{\bf\noindent Problem.}  Fix $T$.
Assume that the set $D$ is  unknown.  Extract information about the
location and shape of $D$ from $u_f(x,t)$, which corresponds to a suitable known
$f$ and is given for all $x\in\partial\Omega$ and $t\in\,]0,\,T[$.

Using the {\it time domain enclosure method} originary developed in \cite{I4}, 
in the previous papers \cite{IK1, IFR, IK2, IK3}, we have considered the problem
above (see also \cite{II} for a system).  The prescribed heat flux $f$ takes the form
$$\begin{array}{ll}
\displaystyle
f(x,t)=\frac{\partial v}{\partial\nu}(x)\varphi(t), & (x,t)\in\partial\Omega\times ]0,\,T[
\end{array}
$$
where $v$ is a special solution of the modified Helmholtz equation $(\Delta-\tau)v=0$, with $\tau>0$, in a domain enclosing
$\overline\Omega$ and, say, $\varphi(t)\sim t^m$ as $t\downarrow 0$, with $m$ being a nonnegative integer.
Note that $f$ depends on $\tau$ and thus, in this sense, the observation data
$u_f(x,t)$, $(x,t)\in\partial\Omega\times\,]0,\,T[$ are infinitely many.  Using the data, we constructed the indicator function
$$\displaystyle
\tau\longmapsto \int_{\partial\Omega}\int_0^Te^{-\tau t}\left(u_f(x,t)\frac{\partial v}{\partial\nu}(x)-v(x)f(x,t)\right)dt\,dS.
$$
From the asymptotic behaviour of this indicator function we extracted several information about the geometry of $D$, more precisely, the distance of an arbitrary point outside of $\Omega$
to $D$, the value of the support function of $D$ at a given direction, the minimum sphere that encloses $D$ with an arbitrary given center point.
Note that in \cite{IFR, IK2}
the governing equation was the heat equation with a variable coefficient
having discontinuity, however, it is easy to see that the present case also
can be covered without difficulty.

Recently, using a new version of the time domain enclosure method developed
in \cite{ESingle}, in \cite{Ithermo}, the author introduced another substitution
of $f$:
$$\begin{array}{ll}
\displaystyle
f(x,t)=\frac{\partial\Theta}{\partial\nu}(x,t), & (x,t)\in\partial\Omega\times\,]0,\,T[,
\end{array}
\tag {1.2}
$$
where $\Theta=\Theta(x,t)$, $(x,t)\in\Bbb R^3\times\,]0,\,T[$, is the solution of the Cauchy problem for the {\it heat equation} 
$$
\displaystyle
(\partial_t-\Delta)\Theta=0
$$
in the whole space, with a special intial data supported on an arbitrary closed ball outside of $\Omega$.  Note that $f$ 
given by (1.2), does not depend on any parameter except for the ball.  The result in \cite{Ithermo} says that the data $u_f$ on $\partial\Omega\times\,]0,\,T[$
yields the distance of the ball $B$ to $D$ from the explicit formula
$$\displaystyle
\lim_{\tau\rightarrow\infty}\frac{1}{\sqrt{\tau}}\log I(\tau)=-2\text{dist}\,(D,B),
$$
where
$$\displaystyle
I(\tau)=\int_{\partial\Omega}
\left(\int_0^Te^{-\tau t}(u_f(x,t)-\Theta (x,t))dt\,\frac{\partial}{\partial\nu}\int_0^{T}e^{-\tau t}\Theta (x,t)dt\right)\,dS
$$
and $\text{dist}\,(D,B)=\inf_{x\in D,\,y\in B}\vert x-y\vert$.

Needless to say, there should be other possibilities of choosing a suitable heat flux which yields the geometry
of $D$ since the function space $L^2(0,T;H^{-1/2}(\partial\Omega))$ is {\it large}.

 Here we have a {\it naive} question:  if one replaces $\Theta$ on (1.2) with a solution of {\it another type of equation}, what happens
on the asymtotic behaviour of the indicator function above?
This is the subject of this paper.  In this paper, as another type of equation we choose the {\it wave equation}.  This choice
comes from the finite propagation speed of the signal governed by the wave equation unlike the heat equation.  How does the wave interact
on the surface of a heat conductive body?  Can one extract information about the geometry of $D$ from the data $u_f$ on $\partial\Omega\times\,]0,\,T[$
by using the idea of the enclosure method under such a choice of $f$?  It should be pointed that in \cite{IK2},
the distance of $\partial\Omega$ to $D$ is also given by using the data $u_f$ on $\partial\Omega\times\,]0,\,T[$ for
an arbitrary $f$ having a {\it positive lower bound}.  However, the $f$ is not necessary a solution of any equation.
See also \cite{IK3} for a result in the enclosure method using a single $f$ which is independent of solutions of any equation.

Now let us describe the results in this paper.
Let $B$ be an open ball satisfying $\overline B\cap\overline\Omega=\emptyset$.
We assume that the radius $\eta$ of $B$ is very small.
Let $\chi_B$ denote the characteristic function of $B$.  Let $v=v_{B,\,\lambda}$ be a solution of
$$\displaystyle
\left\{
\begin{array}{ll}
\displaystyle
(\lambda^2\partial_t^2-\Delta)v=0 & \text{in}\,\Bbb R^3\times\,]0,\,T[,\\
\\
\displaystyle
v(x,0)=0 & \text{in}\,\Bbb R^3,\\
\\
\displaystyle
\partial_tv(x,0)=\Psi_B(x) & \text{in}\,\Bbb R^3,
\end{array}
\right.
\tag {1.3}
$$
where $\lambda$ is a positive constant and
$$\begin{array}{ll}
\displaystyle
\Psi_B(x)=(\eta-\vert x-p\vert)\chi_B(x), & x\in\Bbb R^3
\end{array}
$$
with $p$ denoting the center of $B$.
Note that the function $\Psi_B$ belongs to $H^1(\Bbb R^3)$, since
$\nabla\Psi_B(x)=-\frac{x-p}{\vert x-p\vert}\chi_B(x)$ in the sense of distribution.
The solution $v_B$ of (1.3) is constructed by using the theory of $C_0$-semigroups \cite{Y}.
The class where $v_B$ belongs to is the following:
$$\displaystyle
C^2([0,\,T], L^2(\Bbb R^3))\cap C^1([0,\,T], H^1(\Bbb R^3))\cap C([0,\,T], H^2(\Bbb R^3)).
$$
Needless to say, $v_B$ has an explicit analytical expression, however, we never make use
such expression in the time domain.  We need just the exsitence of $v_B$ in the function spaces indicated above.

The following function is the special $f$ in the problem mentioned above:
$$\begin{array}{ll}
\displaystyle
f_{B,\,\lambda}=f_{B,\,\lambda}(\,\cdot\,,t)=\frac{\partial}{\partial\nu}v_{B,\,\lambda}(\,\cdot\,,t), & t\in\,[0,\,T].
\end{array}
\tag {1.4}
$$

Now we construct the solution $u=u_f$ of (1.1) by prescribing $f=f_{B,\,\lambda}$
and define
$$\begin{array}{lll}
\displaystyle
w_{B,\lambda}(x)=w_{B,\lambda}(x,\tau)=\int_0^Te^{-\tau t}u_f(x,t)\,dt,
& 
x\in\Omega\setminus\overline D,
&
\tau>0,
\end{array}
\tag {1.5}
$$
and
$$\begin{array}{lll}
\displaystyle
w_{B,\,\lambda}^0(x)=w_{B,\,\lambda}^0(x,\tau)=\int_0^Te^{-\tau t}v_{B,\,\lambda}(x,t)\,dt,
&
x\in\Bbb R^3,
&
\tau>0.
\end{array}
\tag {1.6}
$$

We define
$$
\begin{array}{ll}
\displaystyle
I_{\partial\Omega}(\tau;B,\lambda)=\int_{\partial\Omega}
(w_{B,\,\lambda}-w_{B,\,\lambda}^0)\frac{\partial w_{B,\,\lambda}^0}{\partial\nu}\,dS, & \tau>0.
\end{array}
$$
This is the indicator function in the {\it enclosure method} discussed in this paper.
Since we have
$$\displaystyle
\frac{\partial w_{B,\,\lambda}}{\partial\nu}
=\frac{\partial w_{B,\,\lambda}^0}{\partial\nu},
$$
this indicator function has the form
$$\displaystyle
I_{\partial\Omega}(\tau;B,\lambda)=\int_{\partial\Omega}
\left(w_{B,\,\lambda}\frac{\partial w_{B,\,\lambda}^0}{\partial\nu}
-w_{B,\,\lambda}^0\frac{\partial w_{B,\lambda}}{\partial\nu}\right)\,dS.
$$

This indicator function can be computed from the responce $u_f$ on $\partial\Omega$ over
the time interval $]0,\,T[$ which is the solution of (1.1) with $f=f_{B,\,\lambda}$.

\proclaim{\noindent Theorem 1.1(Hardening).}

(i)  If $T$ satisfies
$$\displaystyle
T>\lambda\,\text{dist}(\Omega,\,B)
\tag {1.7}
$$
then there exists a positive number $\tau_0$ such that
$I_{\partial\Omega}(\tau;B,\lambda)>0$ for all $\tau\ge\tau_0$, and we have
$$\displaystyle
\lim_{\tau\longrightarrow\infty}
\frac{1}{\tau}
\log I_{\partial\Omega}(\tau;B,\lambda)=-2\lambda\,\text{dist}\,(\Omega,B).
\tag {1.8}
$$

(ii)  We have
$$\displaystyle
\lim_{\tau\longrightarrow\infty}e^{\tau T}I_{\partial\Omega}(\tau;B,\lambda)=
\left\{
\begin{array}{ll}
\displaystyle
\infty & \text{if}\,\,T>2\lambda\,\text{dist}\,(\Omega,B),\\
\\
\displaystyle
0     & \text{if}\,\,T<2\lambda\,\text{dist}\,(\Omega,B).
\end{array}
\right.
\tag {1.9}
$$

(iii) If $T=2\lambda\,\text{dist}\,(\Omega, B)$, then $e^{\tau T}I_{\partial\Omega}(\tau;B,\lambda)=O(\tau^3)$ as $\tau\longrightarrow\infty$.

\endproclaim

As we already know from \cite{I4,IFR, IK1,IK2, IK3},
there should be various possible choices of the Neumann data $f$ in (1.1) to extract
information about the geometry of an unknown cavity from $u$ on $\partial\Omega\times\,]0,\,T[$.
Theorem 1.1 shows that if one chooses the Neumann data $f=f_{B,\lambda}$ which comes from the solution
of the wave equation with a fixed parameter $\lambda$, then the leading profile of the indicator function does not yield any information
about the cavity.  The choice was bad!  Note also that from (1.8), we have
$$\displaystyle
\lim_{\tau\rightarrow\infty}\frac{1}{\sqrt{\tau}}
\log I_{\partial\Omega}(\tau;B,\lambda)=-\infty.
$$

However, there is another possible choice of $f$ in (1.1).
Given $\tau>0$, choose $\lambda>0$ in (1.4) in such a way that
$$\displaystyle
\lambda^2\tau^2-\tau=0,
\tag {1.10}
$$
that is
$$\displaystyle
\lambda=\frac{1}{\sqrt{\tau}}.
\tag {1.11}
$$
Then, we have the following theorem.

\proclaim{\noindent Theorem 1.2 (Penetrating).}

(i)  Let $T$ be an arbitrary positive number.
Then there exists a positive number $\tau_0$ such that
$I_{\partial\Omega}(\tau;B,\frac{1}{\sqrt{\tau}})>0$ for all $\tau\ge\tau_0$, and we have
$$\displaystyle
\lim_{\tau\longrightarrow\infty}
\frac{1}{\sqrt{\tau}}
\log I_{\partial\Omega}(\tau;B,\frac{1}{\sqrt{\tau}})=-2\,\text{dist}\,(D,B).
\tag {1.12}
$$

(ii)  We have
$$\displaystyle
\lim_{\tau\longrightarrow\infty}e^{\sqrt{\tau}\, T}I_{\partial\Omega}(\tau;B,\frac{1}{\sqrt{\tau}})=
\left\{
\begin{array}{ll}
\displaystyle
\infty & \text{if}\,\,T>2\,\text{dist}\,(D,B),\\
\\
\displaystyle
0     & \text{if}\,\,T<2\,\text{dist}\,(D,B).
\end{array}
\right.
\tag {1.13}
$$

(iii) If $T=2\,\text{dist}\,(D, B)$, then $\displaystyle
e^{\sqrt{\tau}\,T}I_{\partial\Omega}(\tau;B,\frac{1}{\sqrt{\tau}})=O(1)$ as $\tau\longrightarrow\infty$.

\endproclaim

Note that, in this thereom, the Neumann data $f$ is given by (1.4), with $\lambda$ given by (1.11).
Thus, the input data used in Theorem 1.2 vary as $\tau\longrightarrow\infty$ and, in this sense,
they are infinitely many, unlike those of Theorem 1.1.
The first equation on (1.3) with $\lambda$ given by (1.11)
becomes the wave equation with propagation speed $\sqrt{\tau}$:
$$\displaystyle
(\partial_t^2-\tau\Delta)v=0.
\tag {1.14}
$$
So, this sholud be called the enclosure method for the heat equation using a solution of the wave equation with grwoing propagation speed.
One may consider the limit in Theorem 1.2 is a kind of {\it non-relativistic limit}.
Since the speed of (1.14) grows to infinity, we do not need the waiting time for collecting the observation data 
unlike (1.7).  Note that the role of $T$ in (1.12) and (1.13) is different.  In (1.12) $T$ is an arbitrary, however,
to get $\text{dist}\,(D,B)$ by using (1.13) only, we need all $T\in\,]0,\,T_0[$, with $T_0>2\text{dist}\,(D,B)$.

It {\it seems} that equation (1.10) means the vanishing of an {\it obstruction} which prevents the temparture field generated by the flux (1.4)
from entering deep inside of the body $\Omega\setminus\overline D$.  As an evidence, we have different formulae (1.8) and (1.12).
See (2.6) and (2.14) in Section 2 for an explicit role of equation (1.10).

However, we have a question about (1.11).  If $\lambda$ does not satisfy (1.11) exactly, then what happens on the asymptotic behabiour
of the indicator function?  Here instead of (1.11) we choose the case when $\lambda$ is given by
$$\displaystyle
\lambda=\sqrt{\frac{c}{\tau}},
\tag {1.15}
$$
where $c$ is a positive constant.  Then (1.10) becomes
$$\displaystyle
\lambda^2\tau^2-\tau=(c-1)\tau
$$
and equation (1.14)
$$\displaystyle
\left(\partial_t^2-\frac{\tau}{c}\Delta\right)v=0.
$$
This equation also has a grwoing propagation speed as $\tau\longrightarrow\infty$.

In the following result we show that the indicator function $I_{\partial\Omega}(\tau;B,\sqrt{\frac{c}{\tau}}\,)$
has a different asymptotic behaviour across $c=1$.

\proclaim{\noindent Theorem 1.3 (Discontinuity across $c=1$).}

(i)  Let $T$ be an arbitrary positive number.  Let $\pm(c-1)>0$.
Then there exists a positive number $\tau_0$ such that
$\pm I_{\partial\Omega}(\tau;B,\sqrt{\frac{c}{\tau}}\,)>0$ for all $\tau\ge\tau_0$, and we have
$$\displaystyle
\lim_{\tau\longrightarrow\infty}
\frac{1}{\sqrt{\tau}}
\log\left\vert I_{\partial\Omega}(\tau;B,\sqrt{\frac{c}{\tau}}\,)\right\vert
=-2\sqrt{c}\,\text{dist}\,(\Omega,B).
$$

(ii)  Let $\pm(c-1)>0$.  We have
$$\displaystyle
\lim_{\tau\longrightarrow\infty}e^{\sqrt{\tau}\, T}I_{\partial\Omega}(\tau;B,\sqrt{\frac{c}{\tau}}\,)=
\left\{
\begin{array}{ll}
\displaystyle
\pm\infty & \text{if}\,\,T>2\sqrt{c}\,\text{dist}\,(\Omega,B),\\
\\
\displaystyle
0     & \text{if}\,\,T<2\,\sqrt{c}\,\text{dist}\,(\Omega,B).
\end{array}
\right.
$$

(iii) If $T=2\,\sqrt{c}\text{dist}\,(D, B)$, then $\displaystyle
e^{\sqrt{\tau}\,T}I_{\partial\Omega}(\tau;B,\sqrt{\frac{c}{\tau}}\,)=O(1)$ as $\tau\longrightarrow\infty$.

(iv)  Let $c>0$.  We have
$$\displaystyle
c-1\le 
\liminf_{\tau\rightarrow\infty}
\frac{\displaystyle
I_{\partial\Omega}(\tau;B,\sqrt{\frac{c}{\tau}}\,)}
{\displaystyle
\tau\int_{\Omega}\vert w_0\vert^2\,dx}
\le
\limsup_{\tau\rightarrow\infty}
\frac{\displaystyle
I_{\partial\Omega}(\tau;B,\sqrt{\frac{c}{\tau}}\,)}
{\displaystyle
\tau\int_{\Omega}\vert w_0\vert^2\,dx}
\le (c-1)c.
\tag{1.16}
$$

\endproclaim

Note that in (1.16) the case when $c=1$ is also covered.
From Theorems 1.2 and 1.3 we see that the asymptotic behaviour of the indicator function
$I_{\partial\Omega}(\tau;B,\sqrt{\frac{c}{\tau}}\,)$ has a {\it jump discontinuity} at $c=1$.
Only at $c=1$ one can extract information about the cavity from the leading profile of the indicator function
as $\tau\longrightarrow\infty$.

{\bf\noindent Remark 1.4.}
Instead of $f=f_{B,\lambda}$, given by (1.4), prescribe the $f$ in (1.1) as
$$\displaystyle
f=k f_{B,\,\lambda},
$$
where $k$ is a non zero real constant.
The the new indicator function should be
$$\displaystyle
\tilde{I}_{\partial\Omega}(\tau;B,\lambda)
=\int_{\partial\Omega}
\left(\tilde{w}_{B,\,\lambda}\frac{\partial w_{B,\,\lambda}^0}{\partial\nu}-w_{B,\,\lambda}^0\frac{\partial\tilde{w}_{B,\,\lambda}}{\partial\nu}\right)\,dS,
$$
where the function $\displaystyle\tilde{w}_{B,\lambda}$ is given by (1.5) with $f=kf_{B,\lambda}$.
Since we have
$$\displaystyle
u_{kf_{B,\lambda}}=ku_{f_{B,\lambda}},
$$
one gets $\tilde{w}_{B,\lambda}=kw_{B,\lambda}$.  This yields
$$\displaystyle
\tilde{I}_{\partial\Omega}(\tau;B,\lambda)=kI_{\partial\Omega}(\tau;B,\lambda).
$$
Thus, everthing is reduced to studying the case when $k=1$.

Befor closing the introduction, we describe some estimates on $v_{B,\lambda}$ and $u_f$, with $f=f_{B,\lambda}$ given by (1.4), which are employed 
in Section 2.

The $v_{B,\,\lambda}$ is given by a scaling of the classical wave equation.
More precisely,
let $v_0=v_0(x,s)$ solve
$$\displaystyle
\left\{
\begin{array}{ll}
\displaystyle
(\partial_s^2-\Delta)v=0 & \text{in}\,\Bbb R^3\times\,]0,\,\infty[,\\
\\
\displaystyle
v(x,0)=0 & \text{in}\,\Bbb R^3,\\
\\
\displaystyle
\partial_sv(x,0)=\Psi_B(x) & \text{in}\,\Bbb R^3.
\end{array}
\right.
$$
Then the $v_{B,\,\lambda}$ is given by
$$\begin{array}{ll}
\displaystyle
v_{B,\,\lambda}(x,t)=\lambda\,v_0(x,\frac{1}{\lambda}\,t), & (x,t)\in\Bbb R^3\times\,]0,\,T[.
\end{array}
$$
Thus $f_{B,\,\lambda}$, is given by
$$\begin{array}{ll}
\displaystyle
f_{B,\,\lambda}(x,t)=\lambda\,\frac{\partial}{\partial\nu}v_0(x,\frac{1}{\lambda}\,t),
& (x,t)\in\partial\Omega\times\,]0,\,T[.
\end{array}
$$
Using the Fourier transform of $v_0$ with respect to $x\in\Bbb R^3$, we have, for all $s>0$
$$\displaystyle
\frac{1}{s}\Vert v_0(\,\cdot\,,s)\Vert_{L^2(\Bbb R^3)}
+\Vert\partial_s v_0(\,\cdot\,,s)\Vert_{L^2(\Bbb R^3)}
\le 2\Vert\Psi_B\Vert_{L^2(\Bbb R^3)}\equiv C_B.
$$
Thus, one has
$$\displaystyle
\frac{1}{T}\Vert v_{B,\,\lambda}(\,\cdot\,,T)\Vert_{L^2(\Bbb R^3)}+\Vert\partial_t v_{B,\,\lambda}(\,\cdot\,,T)\Vert_{L^2(\Bbb R^3)}\le C_B.
\tag {1.17}
$$
Moreover, the Fourier transform of $v_0$ with respect to $x\in\Bbb R^3$ yields also
$$\displaystyle
\Vert v_0(\,\cdot\,,s)\Vert_{H^1(\Bbb R^3)}\le C_B\sqrt{s^2+3}.
$$
This, together with the trace theorem, yields
$$\displaystyle
\Vert f_{B,\,\lambda}(\,\cdot\,,t)\Vert_{H^{1/2}(\partial\Omega)}
\le C_{\Omega}C_B\sqrt{t^2+3\lambda^2}
$$
and hence
$$\displaystyle
\Vert f_{B,\,\lambda}\Vert_{L^2(0,\,T;H^{-1/2}(\partial\Omega))}
\le C_{\Omega}C_B\sqrt{T^3+3\lambda^2T}.
$$
Thus one gets
$$\displaystyle
\Vert u_f(\,\cdot\,,T)\Vert_{L^2(\Omega\setminus\overline D)}
\le C\sqrt{T^3+3\lambda^2T},
\tag {1.18}
$$
where $C=C_TC_{\Omega}C_B$.  Note that $C$ is independent of $\lambda$.

\section{Proof of Theorems}

In this section, for simplicity of description we always write
$$
\begin{array}{lll}
\displaystyle
w=w_{B,\,\lambda},&
\displaystyle
w_0=w_{B,\,\lambda}^0,&
\displaystyle
R=w-w_0,
\end{array}
$$
where $w_{B,\,\lambda}$ and $w_{B,\,\lambda}^0$ are given by (1.5) and (1.6), respectively.

\subsection{A decomposition formula of the indicator function}

It follows from (1.1) that $w$ satisfies
$$\left\{
\begin{array}{ll}
\displaystyle
(\Delta-\tau)w=e^{-\tau T}u(x,T) & \text{in}\,\Omega\setminus\overline D,\\
\\
\displaystyle
\frac{\partial w}{\partial\nu}=\frac{\partial w_0}{\partial\nu} & \text{on}\,\partial\Omega,\\
\\
\displaystyle
\frac{\partial w}{\partial\nu}=0 & \text{on}\,\partial D.
\end{array}
\right.
$$
Rewrite this as
$$\left\{
\begin{array}{ll}
\displaystyle
(\Delta-\lambda^2\tau^2)w=e^{-\tau T}F& \text{in}\,\Omega\setminus\overline D,\\
\\
\displaystyle
\frac{\partial w}{\partial\nu}=\frac{\partial w_0}{\partial\nu} & \text{on}\,\partial\Omega,\\
\\
\displaystyle
\frac{\partial w}{\partial\nu}=0 & \text{on}\,\partial D,
\end{array}
\right.
\tag {2.1}
$$
where
$$\begin{array}{ll}
\displaystyle
F=F(x,\tau)=u(x,T)+e^{\tau T}(\tau-\lambda^2\tau^2)w,
& x\in\Omega\setminus\overline D.
\end{array}
\tag {2.2}
$$
It follows from (1.3) that the $w_0$ satisfies
$$
\begin{array}{ll}
\displaystyle
(\Delta-\lambda^2\tau^2)w_0+\lambda^2\Psi_B=e^{-\tau T}\lambda^2F_0 & \text{in}\,\Bbb R^3,
\end{array}
\tag {2.3}
$$
where
$$\begin{array}{ll}
\displaystyle
F_0=F_0(x,\tau)
=\partial_tv_{B,\,\lambda}(x,T)+\tau v_{B,\,\lambda}(x,T)
& \text{in $\Bbb R^3$.}
\end{array}
$$
Note that, from (1.17), we have
$$\displaystyle
\Vert F_0\Vert_{L^2(\Bbb R^3)}\le C(1+\tau),
\tag {2.4}
$$
where $C$ is a positive constant independent of $\lambda$.

Then integration by parts, together with (2.1) and (2.3) in $\Omega$, yields
$$
\displaystyle
\int_{\partial\Omega}
\left(\frac{\partial w_0}{\partial\nu}w-\frac{\partial w}{\partial\nu}w_0\right)dS
=\int_{\partial D}w\frac{\partial w_0}{\partial\nu}dS
+e^{-\tau T}\int_{\Omega\setminus\overline D}
(\lambda^2F_0w-Fw_0)dx,
$$
and hence
$$\displaystyle
I_{\partial\Omega}(\tau;B,\lambda)
=\int_{\partial D}w\frac{\partial w_0}{\partial\nu}dS
+e^{-\tau T}\int_{\Omega\setminus\overline D}
(\lambda^2F_0w-Fw_0)dx.
\tag {2.5}
$$
This is the first representation of the indicator function.  
Next we decompose the first term on the right-hand side of
(2.5).  The result yields the following decomposition formula.

\proclaim{\noindent Proposition 2.1.}
We have
$$
\displaystyle
I_{\partial\Omega}(\tau;B,\lambda)=
\int_{\Omega}(\lambda^2\tau^2-\tau)\vert w_0\vert^2\,dx
+J_h(\tau)+E_h(\tau)+{\cal R}(\tau),
\tag {2.6}
$$
where
$$\displaystyle
J_h(\tau)=\int_D(\vert\nabla w_0\vert^2+\tau\vert w_0\vert^2)dx,
\tag {2.7}
$$
$$\displaystyle
E_h(\tau)=\int_{\Omega\setminus\overline D}(\vert\nabla R\vert^2+\tau\vert R\vert^2)dx
\tag {2.8}
$$
and
$$
\displaystyle
{\cal R}(\tau)
=e^{-\tau T}\left\{\int_D\lambda^2F_0w_0dx
+\int_{\Omega\setminus\overline D}u(x,T)Rdx+\int_{\Omega\setminus\overline D}(\lambda^2F_0-u(x,T))w_0dx\right\}.
\tag {2.9}
$$

\endproclaim

{\it\noindent Proof.}
First we show that
$$\begin{array}{ll}
\displaystyle
I_{\partial\Omega}(\tau;B,\lambda)
&
\displaystyle
=J(\tau)+E(\tau)+\tilde{{\cal R}}(\tau),
\end{array}
\tag {2.10}
$$
where
$$\left\{
\begin{array}{l}
\displaystyle
J(\tau)=\int_D(\vert\nabla w_0\vert^2+\lambda^2\tau^2\vert w_0\vert^2)dx,
\\
\\
\displaystyle
E(\tau)=\int_{\Omega\setminus\overline D}(\vert\nabla R\vert^2+\lambda^2\tau^2\vert R\vert^2)dx
\end{array}
\right.
$$
and
$$
\displaystyle
\tilde{{\cal R}}(\tau)
=e^{-\tau T}\left\{\int_D\lambda^2F_0w_0dx
+\int_{\Omega\setminus\overline D}FRdx+\int_{\Omega\setminus\overline D}(\lambda^2F_0-F)w_0dx\right\}.
\tag {2.11}
$$
The proof of (2.10) is now standard in the enclosure method, however, in the next section
we make use of an equation appearing in the proof.  So, for the reader's convenience,
we present the proof.

Since $\overline B\cap\overline\Omega=\emptyset$, we have that
the $R$ satisfies
$$\left\{
\begin{array}{ll}
\displaystyle
(\Delta-\lambda^2\tau^2)R=e^{-\tau T}(F-\lambda^2F_0) & \text{in}\,\Omega\setminus\overline D,\\
\\
\displaystyle
\frac{\partial R}{\partial\nu}=0 & \text{on}\,\partial\Omega,\\
\\
\displaystyle
\frac{\partial R}{\partial\nu}=-\frac{\partial w_0}{\partial\nu} & \text{on}\,\partial D.
\end{array}
\right.
\tag {2.12}
$$
Then one can wite
$$\displaystyle
\int_{\partial D}w\frac{\partial w_0}{\partial\nu}\,dS
=\int_{\partial D}w_0\frac{\partial w_0}{\partial\nu}\,dS-\int_{\partial D}R\frac{\partial R}{\partial\nu}\,dS.
$$
It follows from (2.3) in $D$ that
$$
\displaystyle
\int_{\partial D}w_0\frac{\partial w_0}{\partial\nu}\,dS
=\int_D(\vert\nabla w_0\vert^2+\lambda^2\tau^2\vert w_0\vert^2)dx
+e^{-\tau T}\int_D\lambda^2F_0w_0dx.
$$
It follows from (2.12) that
$$\begin{array}{ll}
\displaystyle
-\int_{\partial D}R\frac{\partial R}{\partial\nu} dS
&
\displaystyle
=\int_{\partial\,(\Omega\setminus\overline D)}R\frac{\partial R}{\partial\nu} dS\\
\\
\displaystyle
&
\displaystyle
=\int_{\Omega\setminus\overline D}(\vert\nabla R\vert^2+\lambda^2\tau^2\vert R\vert^2)dx
+e^{-\tau T}\int_{\Omega\setminus\overline D}(F-\lambda^2F_0)Rdx.
\end{array}
\tag {2.13}
$$
Thus, we obtain
$$\begin{array}{l}
\displaystyle
\,\,\,\,\,\,
\int_{\partial D}w\frac{\partial w_0}{\partial\nu}\,dS\\
\\
\displaystyle
=\int_D(\vert\nabla w_0\vert^2+\lambda^2\tau^2\vert w_0\vert^2)dx+
\int_{\Omega\setminus\overline D}(\vert\nabla R\vert^2+\lambda^2\tau^2\vert R\vert^2)dx\\
\\
\displaystyle
\,\,\,
+e^{-\tau T}\left\{\int_D\lambda^2F_0w_0dx+\int_{\Omega\setminus\overline D}(F-\lambda^2F_0)Rdx\right\}.
\end{array}
$$
Then a combination of this and (2.5) yields (2.10).

Substituting (2.2) into the right-hand side on (2.11), we can rewrite (2.11) as
$$\begin{array}{ll}
\displaystyle
\,\,\,\,\,\,
\tilde{{\cal R}}(\tau)
&
\displaystyle
=e^{-\tau T}\left\{\int_D\lambda^2F_0w_0dx
+\int_{\Omega\setminus\overline D}u(x,T)Rdx+\int_{\Omega\setminus\overline D}(\lambda^2F_0-u(x,T))w_0dx\right\}
\\
\\
\displaystyle
&
\displaystyle
\,\,\,
+\int_{\Omega\setminus\overline D}(\tau-\lambda^2\tau^2)(\vert R\vert^2-\vert w_0\vert^2)\,dx\\
\\
\displaystyle
&
\displaystyle
=
e^{-\tau T}\left\{\int_D\lambda^2F_0w_0dx
+\int_{\Omega\setminus\overline D}u(x,T)Rdx+\int_{\Omega\setminus\overline D}(\lambda^2F_0-u(x,T))w_0dx\right\}
\\
\\
\displaystyle
&
\displaystyle
\,\,\,
+\int_{\Omega\setminus\overline D}(\tau-\lambda^2\tau^2)\vert R\vert^2\,dx
-\int_{\Omega}(\tau-\lambda^2\tau^2)\vert w_0\vert^2\,dx
+\int_D(\tau-\lambda^2\tau^2)\vert w_0\vert^2\,dx.
\end{array}
$$
Then (2.10) becoms (2.6).

\noindent
$\Box$

\subsection{Estimating indicator functions}

First we give a {\it rough} estimate of $E_h(\tau)$ from above in terms of $J_h(\tau)$ and $\Vert w_0\Vert_{L^2(\Omega\setminus\overline D)}$.
\proclaim{\noindent Lemma 2.2.}  Let $\epsilon>0$.
We have, as $\tau\longrightarrow\infty$
$$\begin{array}{ll}
\displaystyle
E_h(\tau) & 
\displaystyle
\le C_1(\lambda^4\tau^3+1)J_h(\tau)
+C_2(\epsilon)M(\lambda,\tau)
e^{-2\tau T}\\
\\
\displaystyle
&
\,\,\,
\displaystyle
+\frac{1+\epsilon}{\tau}\int_{\Omega\setminus\overline D}(\lambda^2\tau^2-\tau)^2\vert w_0\vert^2\,dx,
\end{array}
\tag {2.14}
$$
where $C_1$ and $C_2(\epsilon)$ are positive constants inedependent of $\lambda$ and $\tau$, and
$$\displaystyle
M(\lambda,\tau)=\lambda^4(1+\tau)^2+\frac{1}{\tau}\{(1+\lambda^2)+\lambda^4(1+\tau)^2\}.
$$
\endproclaim

{\it\noindent Proof.}
It follows from the boundary condition on $\partial D$ in (2.12) and (2.13) that
$$\begin{array}{l}
\displaystyle
\,\,\,\,\,\,
\int_{\Omega\setminus\overline D}\left(\vert\nabla R\vert^2+\lambda^2\tau^2\vert R\vert^2+e^{-\tau T}(F-\lambda^2F_0)R\right)\,dx
\\
\\
\displaystyle
=\int_{\partial D}\frac{\partial w_0}{\partial\nu}RdS.
\end{array}
$$
Using (2.2), we have
$$\begin{array}{l}
\displaystyle
\,\,\,\,\,\,
\lambda^2\tau^2\vert R\vert^2+e^{-\tau T}(F-\lambda^2F_0)R\\
\\
\displaystyle
=\lambda^2\tau^2\vert R\vert^2+(\tau-\lambda^2\tau^2)wR+e^{-\tau T}(u(x,T)-\lambda^2F_0)R
\\
\\
\displaystyle
=\lambda^2\tau^2\vert R\vert^2+(\tau-\lambda^2\tau^2)(w_0+R)R+e^{-\tau T}(u(x,T)-\lambda^2F_0)R
\\
\\
\displaystyle
=\tau\vert R\vert^2+(\tau-\lambda^2\tau^2)w_0 R+e^{-\tau T}(u(x,T)-\lambda^2F_0)R\\
\\
\displaystyle
=\tau\left\vert R+\frac{(\tau-\lambda^2\tau^2)w_0+e^{-\tau T}(u(x,T)-\lambda^2F_0)}
{2\tau}\right\vert^2\\
\\
\displaystyle
\,\,\,
-\frac{\vert(\tau-\lambda^2\tau^2)w_0+e^{-\tau T}(u(x,T)-\lambda^2F_0)\vert^2}{4\tau},
\end{array}
$$
and thus
$$\begin{array}{l}
\displaystyle
\,\,\,\,\,\,
\int_{\Omega\setminus\overline D}
\left(\vert\nabla R\vert^2
+\tau\left\vert R+\frac{(\tau-\lambda^2\tau^2)w_0+e^{-\tau T}(u(x,T)-\lambda^2F_0)}{2\tau}\right\vert^2\right)\,dx\\
\\
\displaystyle
=\int_{\partial D}\frac{\partial w_0}{\partial\nu}RdS
+\frac{1}{4\tau}\int_{\Omega\setminus\overline D}\vert(\tau-\lambda^2\tau^2)w_0+e^{-\tau T}(u(x,T)-\lambda^2F_0)\vert^2\,dx.
\end{array}
$$
This yields
$$\begin{array}{l}
\displaystyle
\,\,\,\,\,\,
\frac{1}{2}\int_{\Omega\setminus\overline D}
(\vert\nabla R\vert^2+\tau\vert R\vert^2)\,dx
\\
\\
\displaystyle
\le
\int_{\partial D}\frac{\partial w_0}{\partial\nu}\,R\,dS
+\frac{1}{2\tau}\int_{\Omega\setminus\overline D}\vert(\tau-\lambda^2\tau^2)w_0+e^{-\tau T}(u(x,T)-\lambda^2F_0)\vert^2\,dx,
\end{array}
$$
and thus
$$\begin{array}{l}
\displaystyle
\,\,\,\,\,\,
\int_{\Omega\setminus\overline D}
(\vert\nabla R\vert^2+\tau\vert R\vert^2)\,dx\\
\\
\displaystyle
\le
2\int_{\partial D}\frac{\partial w_0}{\partial\nu}\,R\,dS
+\frac{1}{\tau}\int_{\Omega\setminus\overline D}\vert(\tau-\lambda^2\tau^2)w_0+e^{-\tau T}(u(x,T)-\lambda^2F_0)\vert^2\,dx.
\end{array}
$$
From (1.18) and (2.4) we have 
$$
\displaystyle
\Vert u(x,T)-\lambda^2F_0\Vert_{L^2(\Omega\setminus\overline D)}\le C_3(\sqrt{1+\lambda^2}+\lambda^2(1+\tau)),
$$
where $C_3$ is a positive constant independent of $\lambda$ and $\tau$.
This, together with (2.8) and the inequality
$$\begin{array}{lll}
\displaystyle
(a+b)^2\le (1+\epsilon)a^2+(1+4\epsilon^{-1})b^2, & a>0, & b>0,
\end{array}
$$
yields
$$\begin{array}{l}
\displaystyle
\,\,\,\,\,\,
E_h(\tau)
\le 2\int_{\partial D}\frac{\partial w_0}{\partial\nu}RdS+
\frac{1+\epsilon}{\tau}\int_{\Omega\setminus\overline D}(\lambda^2\tau^2-\tau)^2\vert w_0\vert^2\,dx
+C_{\epsilon}C_3^2K(\lambda,\tau)e^{-2\tau T},
\end{array}
\tag {2.15}
$$
where $C_{\epsilon}=1+4\epsilon^{-1}$ and
$$\displaystyle
K(\lambda,\tau)=\frac{1}{\tau}(\sqrt{1+\lambda^2}+\lambda^2(1+\tau))^2.
$$
By the trace theorem \cite{Gr}, one can choose a positive constant $C=C(D,\Omega)$ and $\tilde{R}\in H^1(D)$
such that $\tilde{R}=R$ on $\partial D$ and $\Vert\tilde{R}\Vert_{H^1(D)}\le C\Vert R\Vert_{H^1(\Omega\setminus\overline D)}$.
Then, we have
$$\begin{array}{l}
\displaystyle
\,\,\,\,\,\,
\int_{\partial D}\frac{\partial w_0}{\partial\nu}RdS\\
\\
\displaystyle
=\int_{\partial D}\frac{\partial w_0}{\partial\nu}\tilde{R}dS\\
\\
\displaystyle
=\int_D(\Delta w_0)\tilde{R}\,dx+\int_D\nabla w_0\cdot\nabla\tilde{R}\,dx\\
\\
\displaystyle
=\int_D\lambda^2\tau^2w_0\tilde{R}dx+\int_D\nabla w_0\cdot\nabla\tilde{R}\,dx+e^{-\tau T}\int_D\lambda^2F_0\tilde{R}dx.
\end{array}
$$
Note that in the last step, we have made use of equation (2.3) on $D$.
Then the choice of $\tilde{R}$ and (2.4) yield
$$\begin{array}{l}
\,\,\,\,\,\,
\displaystyle
\left\vert\int_{\partial D}\frac{\partial w_0}{\partial\nu}\,R\,dS
\right\vert\\
\\
\displaystyle
\le
C\Vert R\Vert_{H^1(\Omega\setminus\overline D)}
\left(\lambda^2\tau^2\Vert w_0\Vert_{L^2(D)}+\Vert\nabla w_0\Vert_{L^2(D)}+C_4\lambda^2(1+\tau)e^{-\tau T}\right).
\end{array}
\tag {2.16}
$$
Here we note that $\Vert R\Vert_{H^1(\Omega\setminus\overline D)}\le E_h(\tau)^{1/2}$ for all $\tau\ge 1$, and
$\Vert w_0\Vert_{L^2(D)}\le\tau^{-1/2}J_h(\tau)^{1/2}$, $\Vert\nabla w_0\Vert_{L^2(D)}\le J_h(\tau)^{1/2}$ for all $\tau>0$.
From these, (2.15) and (2.16), we obtain
$$\begin{array}{ll}
\displaystyle
E_h(\tau)
&
\displaystyle
\le
C'E_h(\tau)^{1/2}
\left\{(\lambda^2\tau^{3/2}+1)J_h(\tau)^{1/2}+C_4\lambda^2(1+\tau) e^{-\tau T}\right\}\\
\\
\displaystyle
\,\,\,
&
\displaystyle
\,\,\,
+\frac{1+\epsilon}{\tau}
\int_{\Omega\setminus\overline D}(\lambda^2\tau^2-\tau)^2\vert w_0\vert^2\,dx
+
C_{\epsilon}C_{3}^2K(\lambda,\tau)e^{-2\tau T},
\end{array}
$$
where $C'$ is a positive constant independent of $\lambda$ and $\tau$.
Now a standard argument yields (2.14).

\noindent
$\Box$

{\bf\noindent Remark 2.3.}  It is very important to have the factor $1+\epsilon$ in the third term of the right-hand side on (2.14).
This yields the upper bound (2.30) for the indicator function $I_{\partial\Omega}(\tau;B,\sqrt{\frac{c}{\tau}}\,)$.

Next we describe local upper and lower estimates for $w_0$.

\proclaim{\noindent Lemma 2.4 (Propagation estimates).} Let $U$ be an arbitrary bounded open subset of $\Bbb R^3$ such that $\overline B\cap\overline U=\emptyset$.

(i)  We have
$$
\displaystyle
\lambda\tau\Vert w_0\Vert_{L^2(U)}+
\Vert\nabla w_0\Vert_{L^2(U)}
\le C\left\{\lambda^3\tau e^{-\tau\lambda\,\text{dist}\,(U,B)}+\frac{\lambda(1+\tau)}{\tau}\,e^{-\tau T}\right\},
\tag {2.17}
$$
where $C$ is a positive constant independent of $\lambda$ and $\tau$.

(ii) Let $\partial U$ be $C^2$.
Fix $\lambda$ and let $T$ satisfy
$$\displaystyle
T>\lambda\text{dist}\,(U,B).
\tag {2.18}
$$
Then there exist positive constants $\tau_0$ and $C$ such that, for all $\tau\ge\tau_0$
$$\displaystyle
\tau^{12}e^{2\tau\lambda\,\text{dist}\,(U,B)}\int_{U}\vert w_0\vert^2\,dx\ge C.
\tag {2.19}
$$

(iii)  Let $\partial U$ be $C^2$.
Let $\displaystyle\lambda=\sqrt{\frac{c}{\tau}}$ with a positive constant $c$.
Let $T$ be an arbitraly positive number.
Then there exist positive constants $\tau_0$ and $C$ such that, for all $\tau\ge\tau_0$
$$\displaystyle
\tau^{8}e^{2\sqrt{\tau}\,\sqrt{c}\,\text{dist}\,(U,B)}\int_{U}\vert w_0\vert^2\,dx\ge C.
\tag {2.20}
$$

\endproclaim

{\it\noindent Proof.}
We set
$$
\displaystyle
\epsilon_0=e^{\tau T}(w_0-v_0),
$$
where $v_0\in H^1(\Bbb R^3)$ is the solution of
$$\begin{array}{ll}
\displaystyle
(\Delta-\lambda^2\tau^2)v_0+\lambda^2\Psi_B=0
&
\text{in $\Bbb R^3$.}
\end{array}
$$
The $v_0$ has the explicit form
$$
\displaystyle
v_0(x)
=
\frac{\lambda^2}{4\pi}
\int_B
\frac{e^{-\tau\lambda\vert x-y\vert}}{\vert x-y\vert}(\eta-\vert y-p\vert)\,dy.
\tag {2.21}
$$
We have
$$\displaystyle
w_0=v_0+e^{-\tau T}\epsilon_0
$$
and, from (2.3),
$$\begin{array}{ll}
\displaystyle
(\Delta-\lambda^2\tau^2)\epsilon_0=\lambda^2F_0 & \text{in}\,\Bbb R^3.
\end{array}
\tag {2.22}
$$
Then, from (2.4) and (2.22), we can easily see that 
$$
\displaystyle 
\lambda\tau\Vert\epsilon_0\Vert_{L^2(\Bbb R^3)}+\Vert\nabla\epsilon_0\Vert_{L^2(\Bbb R^3)}\le C_5\frac{\lambda(1+\tau)}{\tau},
\tag {2.23}
$$
where $C_5$ is a positive constant independent of $\lambda$ and $\tau$.

The expression (2.21) for $v_0$ yields
$$
\displaystyle
\lambda\tau\Vert v_0\Vert_{L^2(U)}
+
\Vert\nabla v_0\Vert_{L^2(U)}
\le C_6\lambda^3\tau e^{-\tau\lambda\,\text{dist}\,(U,B)},
\tag {2.24}
$$
where $C_6$ is a positive constant independent of $\lambda$ and $\tau$.
A combination of (2.23) and (2.24) gives (2.17).

It follows from (2.23) that
$$
\displaystyle
\int_{U}\vert w_0\vert^2\,dx\ge
\frac{1}{2}\int_{U}\vert v_0\vert^2\,dx-C_5^2\left(\frac{1+\tau}{\tau}\right)^2\frac{e^{-2\tau T}}{\tau^2}.
$$
From Appendix in \cite{ESingle}, we know that (2.21) has the expression
$$
\displaystyle
v_0(x)
=\frac{1}{\tau^4\lambda^2}\frac{e^{-\tau\lambda\vert x-p\vert}}{\vert x-p\vert}
(-2\cosh(\tau\lambda\eta)+\tau\lambda\eta\sinh(\tau\lambda\eta)+2).
$$
This yields, for sufficiently large $\tau$,
$$\displaystyle
\int_{U}\vert v_0\vert^2\,dx
\ge
\left\{\begin{array}{ll}
\displaystyle
C_7^2\tau^{-6}\int_{U}\frac{e^{-2\tau\lambda\,(\vert x-p\vert-\eta)}}{\vert x-p\vert^2}\,dx
&
\text{if $\lambda$ is fixed,}
\\
\\
\displaystyle
C_8^2\tau^{-5}\int_{U}\frac{e^{-2\sqrt{\tau}\,\sqrt{c}\,(\vert x-p\vert-\eta)}}{\vert x-p\vert^2}\,dx
&
\text{if $\displaystyle\lambda=\sqrt{\frac{c}{\tau}}$,}
\end{array}
\right.
\tag {2.25}
$$
where $C_7$ is a positive constant independent of $\tau$, and $C_8$ is independent of $\tau$ and $\lambda$.
Assume that $\partial U$ is $C^2$.
In \cite{IW00, IEE}, we have already proved that there exist positive constants $\tau_0$ and $C'$ such that
for all $\tau\ge\tau_0$ 
$$\displaystyle
\tau^6 e^{2\tau\,\text{dist}\,(U,B)}\int_{U}\frac{e^{-2\tau(\vert x-p\vert-\eta)}}{\vert x-p\vert^2}\,dx\ge C',
$$
and thus
$$\displaystyle
\tau^3 e^{2\sqrt{\tau}\,\sqrt{c}\,\text{dist}\,(U,B)}\int_{U}\frac{e^{-2\sqrt{\tau}\,\sqrt{c}\,(\vert x-p\vert-\eta)}}{\vert x-p\vert^2}\,dx\ge C'.
$$
Now applying these to the right-hand side on (2.25), 
we see the validity of (ii) and (iii).

\noindent
$\Box$

{\bf\noindent Remark 2.5.}
In the proof of (2.17), the estimate (2.24) is essential.  For this purpose, we made use of 
the expression of $v_0$ given by (2.21) only.

Here we describe preliminary estimates for the indicator function.

\proclaim{\noindent Lemma 2.6.} Let $T$ be an arbitrary positive number.

(i) Fix $\lambda>0$.
We have, as $\tau\longrightarrow\infty$
$$
\displaystyle
I_{\partial\Omega}(\tau;B,\lambda)=
O(\tau^3e^{-2\tau\lambda\,\text{dist}\,(\Omega,B)}+\tau^3 e^{-2\tau T}
+\tau e^{-\tau T}e^{-\tau\lambda\,\text{dist}\,(\Omega,B)})
\tag {2.26}
$$
and
$$
\displaystyle
I_{\partial\Omega}(\tau;B,\lambda)
\ge \int_{\Omega}(\lambda^2\tau^2-\tau)\vert w_0\vert^2\,dx
+O(\tau e^{-\tau T}e^{-\tau\lambda\,\text{dist}\,(\Omega,B)}
+\tau e^{-2\tau T}).
\tag {2.27}
$$

(ii) 
We have, as $\tau\longrightarrow\infty$
$$\displaystyle
I_{\partial\Omega}(\tau;B,\frac{1}{\sqrt{\tau}})
=O(e^{-2\sqrt{\tau}\,\text{dist}\,(D,B)}+\tau^{-1}e^{-\tau T}e^{-\sqrt{\tau}\,\text{dist}\,(\Omega,B)})
\tag {2.28}
$$
and
$$\displaystyle
I_{\partial\Omega}(\tau;B,\frac{1}{\sqrt{\tau}})
\ge \tau\int_D\vert w_0\vert^2\,dx+O(\tau^{-1}e^{-\tau T}e^{-\sqrt{\tau}\,\text{dist}\,(\Omega,B)}).
\tag {2.29}
$$

(iii)  Let $\epsilon$ be an arbitrary positive number.  We have, as $\tau\longrightarrow\infty$
$$\begin{array}{l}
\displaystyle
\,\,\,\,\,\,
I_{\partial\Omega}(\tau;B,\sqrt{\frac{c}{\tau}}\,)\\
\\
\displaystyle
\le
\tau (1+\epsilon)(c-1)\left(c-\frac{\epsilon}{1+\epsilon}\right)\int_{\Omega}\vert w_0\vert^2dx
+O(\tau e^{-\tau T}e^{-2\sqrt{\tau}\,\sqrt{c}\,\text{dist}\,(\Omega,B)})
\end{array}
\tag {2.30}
$$
and
$$\begin{array}{l}
\displaystyle
\,\,\,\,\,\,
I_{\partial\Omega}(\tau;B,\sqrt{\frac{c}{\tau}}\,)\\
\\
\displaystyle
\ge
\tau (c-1)\int_{\Omega}\vert w_0\vert^2dx
+O(\tau e^{-\tau T}e^{-2\sqrt{\tau}\,\sqrt{c}\,\text{dist}\,(\Omega,B)}).
\end{array}
\tag {2.31}
$$

\endproclaim

{\it\noindent Proof.}
First we give a proof of (2.26) and (2.27).
From (2.17), with $\star=D,\Omega\setminus\overline D$ we have
$$
\displaystyle
\Vert w_0\Vert_{L^2(\star)}
=O(e^{-\tau\lambda\,\text{dist}\,(\star,B)}+\tau^{-1}e^{-\tau T})
\tag {2.32}
$$
and
$$
\displaystyle
\Vert\nabla w_0\Vert_{L^2(\star)}
=O(\tau e^{-\tau\lambda\,\text{dist}\,(\star,B)}+e^{-\tau T}).
$$
Thus, (2.7) gives
$$
\displaystyle
J_h(\tau)=O(\tau^2 e^{-2\tau\lambda\,\text{dist}\,(D,B)}
+e^{-2\tau T}).
\tag {2.33}
$$
Then, from (2.14), (2.32) with $\star=\Omega\setminus\overline D$ and (2.33), we have
$$
\displaystyle
E_h(\tau)
=O(\tau^5e^{-2\tau\lambda\,\text{dist}\,(D,B)}
+\tau^3e^{-2\tau\lambda\,\text{dist}\,(\Omega,B)}+\tau^3 e^{-2\tau T}).
$$
Since $\text{dist}\,(D,B)>\text{dist}\,(\Omega,B)$, this yields
$$
\displaystyle
E_h(\tau)
=O(\tau^3e^{-2\tau\lambda\,\text{dist}\,(\Omega,B)}+\tau^3 e^{-2\tau T}).
\tag {2.34}
$$
Since $\tau\Vert R\Vert_{L^2(\Omega\setminus\overline D)}^2\le E_h(\tau)$, 
(2.34) gives
$$\displaystyle
\Vert R\Vert_{L^2(\Omega\setminus\overline D)}
=O(\tau e^{-\tau\lambda\,\text{dist}\,(\Omega,B)}+\tau e^{-\tau T}).
$$
This, together with (1.18), gives
$$\displaystyle
\int_{\Omega\setminus\overline D}u(x,T)Rdx
=O(\tau e^{-\tau\lambda\,\text{dist}\,(\Omega,B)}+\tau e^{-\tau T}).
$$
From (2.32), (1.18) and (2.4), we obtain
$$
\displaystyle
\int_D\lambda^2F_0w_0dx=O(\tau e^{-\tau\lambda\,\text{dist}\,(D,B)}+e^{-\tau T})
$$
and
$$
\displaystyle
\int_{\Omega\setminus\overline D}(\lambda^2F_0-u(x,T))w_0dx
=O(\tau e^{-\tau\lambda\,\text{dist}\,(\Omega,B)}+e^{-\tau T}).
$$
Applying these to the right-hand side on (2.9), we obtain
$$
\displaystyle
{\cal R}(\tau)
=O(\tau e^{-\tau T}e^{-\tau\lambda\,\text{dist}\,(\Omega,B)}
+\tau e^{-2\tau T}).
\tag {2.35}
$$
Now applying (2.33), (2.34) and (2.35) to (2.6), we obtain 
$$\begin{array}{l}
\displaystyle
\,\,\,\,\,\,
I_{\partial\Omega}(\tau;B,\lambda)\\
\\
\displaystyle
=\int_{\Omega}(\lambda^2\tau^2-\tau)\vert w_0\vert^2\,dx\\
\\
\displaystyle
\,\,\,
+O(\tau^3e^{-2\tau\lambda\,\text{dist}\,(\Omega,B)}+\tau^3 e^{-2\tau T}
+\tau e^{-\tau T}e^{-\tau\lambda\,\text{dist}\,(\Omega,B)}),
\end{array}
$$
and this and (2.32) with $\star=\Omega$ yield (2.26).

By omitting $J_h(\tau)$ and $E_h(\tau)$ in (2.6) which are non negative, we have
$$\displaystyle
I_{\partial\Omega}(\tau;B,\lambda)
\ge
\int_{\Omega}(\lambda^2\tau^2-\tau)\vert w_0\vert^2\,dx+{\cal R}(\tau).
$$
Now from (2.35) we obtain (2.27).

Next we give a proof of (2.28), (2.29), (2.30) and (2.31).
Let $\displaystyle\lambda$ be given by (1.15), which covers (1.11) as a special case.
From (2.17) with $\star=D,\Omega\setminus\overline D$ we have
$$
\displaystyle
\Vert w_0\Vert_{L^2(\star)}
=O(\tau^{-1}e^{-\sqrt{\tau}\,\sqrt{c}\,\text{dist}\,(\star,B)}+\tau^{-1}e^{-\tau T})
\tag {2.36}
$$
and
$$
\displaystyle
\Vert\nabla w_0\Vert_{L^2(\star)}
=O(\frac{1}{\sqrt{\tau}}\tau e^{-\sqrt{\tau}\,\sqrt{c}\,\text{dist}\,(\star,B)}+\frac{1}{\sqrt{\tau}}e^{-\tau T}).
$$
Thus, one gets
$$
\displaystyle
J_h(\tau)=O(\tau^{-1}e^{-2\sqrt{\tau}\,\sqrt{c}\,\text{dist}\,(D,B)}
+\tau^{-1}\,e^{-2\tau T}).
\tag {2.37}
$$
Moreover, from (2.14) we have
$$\displaystyle
E_h(\tau)\le C_1(\tau+1)J_h(\tau)+
(1+\epsilon)(c-1)^2\tau\int_{\Omega}\vert w_0\vert^2\,dx+
O(e^{-2\tau T}),
\tag {2.38}
$$
and thus 
$$\displaystyle
E_h(\tau)=O(e^{-2\sqrt{\tau}\,\text{dist}\,(D,B)}+e^{-2\tau T}+(c-1)^2\tau^{-1}e^{-2\sqrt{\tau}\,\sqrt{c}\,\text{dist}\,(\Omega,B)}).
\tag {2.39}
$$
This yields
$$\displaystyle
\Vert R\Vert_{L^2(\Omega\setminus\overline{D})}=
O(\frac{e^{-\sqrt{\tau}\,\text{dist}\,(D,B)}}{\sqrt{\tau}}+\frac{e^{-\tau T}}{\sqrt{\tau}}+
\vert c-1\vert\tau^{-1}e^{-\sqrt{\tau}\,\sqrt{c}\,\text{dist}\,(\Omega,B)}).
$$
Then, from (2.9) we have
$$\begin{array}{l}
\displaystyle
\,\,\,\,\,\,
\cal R(\tau)\\
\\
\displaystyle
=O(e^{-\tau T}\left(\frac{1}{\tau}e^{-\sqrt{\tau}\,\sqrt{c}\,\text{dist}\,(D,B)}
+\frac{1}{\tau}e^{-\tau T}+\frac{e^{-\sqrt{\tau}\,\sqrt{c}\,\text{dist}\,(D,B)}}{\sqrt{\tau}}+\frac{e^{-\tau T}}{\sqrt{\tau}}
\right.
\\
\\
\displaystyle
\,\,\,
\left.
+\frac{1}{\tau}e^{-\sqrt{\tau}\,\sqrt{c}\,\text{dist}\,(\Omega,B)}
+\frac{1}{\tau}e^{-\tau T}\,\right))\\
\\
\displaystyle
\,\,\,
+O(\vert c-1\vert\tau^{-1}e^{-\tau T}e^{-\sqrt{\tau}\,\sqrt{c}\,\text{dist}\,(\Omega,B)})
\\
\\
\displaystyle
=O(\tau^{-1}e^{-\tau T}e^{-\sqrt{\tau}\,\sqrt{c}\,\text{dist}\,(\Omega,B)})
+O(\vert c-1\vert\tau^{-1}e^{-\tau T}e^{-\sqrt{\tau}\,\sqrt{c}\,\text{dist}\,(\Omega,B)})\\
\\
\displaystyle
=O(\tau^{-1}e^{-\tau T}e^{-\sqrt{\tau}\,\sqrt{c}\,\text{dist}\,(\Omega,B)}).
\end{array}
\tag {2.40}
$$
From (2.6) and (2.40), we have
$$\displaystyle
I_{\partial\Omega}(\tau;B,\sqrt{\frac{c}{\tau}}\,)
=J_h(\tau)+\tau(c-1)\int_{\Omega}\vert w_0\vert^2\,dx+E_h(\tau)+O(\tau^{-1}e^{-\tau T}e^{-\sqrt{\tau}\,\sqrt{c}\,\text{dist}\,(\Omega,B)}).
\tag {2.41}
$$
Let $c=1$.  Appllying (2.36) with $\star=\Omega$, (2.37) and (2.39) to the right-hand side on (2.41), we obtain (2.28).  
(2.29) is now clear.

Next consider the case when $c\not=1$. 
Since we have
$$\displaystyle
(c-1)+(1+\epsilon)(c-1)^2
=(1+\epsilon)(c-1)\left(c-\frac{\epsilon}{1+\epsilon}\right),
$$ 
from (2.38), one gets
$$\begin{array}{l}
\displaystyle
\,\,\,\,\,\,
\tau(c-1)\int_{\Omega}\vert w_0\vert^2\,dx+E_h(\tau)\\
\\
\displaystyle
\le C_1(\tau+1)J_h(\tau)+
(1+\epsilon)(c-1)\left(c-\frac{\epsilon}{1+\epsilon}\right)\,\tau\int_{\Omega}\vert w_0\vert^2dx
+O(e^{-2\tau T}).
\end{array}
$$
This, together with (2.41) and (2.37), yields (2.30). 
(2.31) is an easy consequence of (2.41).

\noindent
$\Box$

\subsection{Proof of Theorems 1.1 and 1.2}

First we give a proof of Theorem 1.1.
Let $T\le 2\lambda\,\text{dist}\,(\Omega,\,B)$.
It follows from (2.26) that
$$
\displaystyle
e^{\tau T}I_{\partial\Omega}(\tau;B,\lambda)
=O(\tau^3e^{-\tau(2\lambda\,\text{dist}\,(\Omega,B)-T)}+\tau^3 e^{-\tau T}
+\tau e^{-\tau\lambda\,\text{dist}\,(\Omega,B)}).
$$
This yields (1.9) in the case of $T<2\lambda\,\text{dist}\,(\Omega,B)$
and (iii).

Next let $T$ satisfy (1.7), i.e., (2.18) with $U=\Omega$.
Then, we have (2.19) with $U=\Omega$.
This, together with (2.27), yields that there exist positive constants $C$ and $\tau_0$ such that for all $\tau\ge\tau_0$,
$$\displaystyle
\tau^{10}e^{2\tau\,\lambda\,\text{dist}\,(\Omega,B)}I_{\partial\Omega}(\tau;B,\lambda)\ge
C+O(\tau^{11} e^{-\tau(T-\lambda\,\text{dist}\,(\Omega,B))}).
\tag {2.42}
$$
A combination of (2.26) and (2.42) yields (1.8).  (ii) in the case of $T>2\lambda\,\text{dist}\,(\Omega,B)$
is a direct consequence of (1.8), since $T$ satisfies (1.7) in this case.

The proof of Theorem 1.2 is as follows.
From (2.20) with $U=D$ and (2.29), we have
$$\displaystyle
\tau^7 e^{2\sqrt{\tau}\,\text{dist}\,(D,B)}I_{\partial\Omega}(\tau;B,\frac{1}{\sqrt{\tau}})
\ge C+O(\tau^{6}e^{-\tau T}e^{\sqrt{\tau}\,(\text{dist}\,(D,B)-\text{dist}\,(\Omega,B))}\,).
\tag {2.43}
$$
Note that there is no retsriction on $T$.
Now it is easy to see that a combination of (2.28) and (2.43) yields the validity of Theorem 1.2.

\subsection{Proof of Theorem 1.3}  Using (2.30) and (2.31), we see that
(i), (ii) and (iii) can be easily derived as those of Theorem 1.2.  
(1.16) in the case when $c=1$ is a direct consequence of Theorem 1.2 (i) and (2.20) with $U=\Omega$.
Let $c\not=1$.  It follows from (2.20), with $U=\Omega$, (2.30) and (2.31) that
$$
\displaystyle
\limsup_{\tau\rightarrow\infty}
\frac{\displaystyle
I_{\partial\Omega}(\tau;B,\sqrt{\frac{c}{\tau}}\,)}
{\displaystyle
\tau\int_{\Omega}\vert w_0\vert^2\,dx}\le (1+\epsilon)(c-1)
\left(c-\frac{\epsilon}{1+\epsilon}\right)
\tag {2.44}
$$
and
$$\displaystyle
\liminf_{\tau\rightarrow\infty}
\frac{\displaystyle
I_{\partial\Omega}(\tau;B,\sqrt{\frac{c}{\tau}}\,)}
{\displaystyle
\tau\int_{\Omega}\vert w_0\vert^2\,dx}\ge c-1.
$$
Since the left hand-side on (2.44) is independent of $\epsilon$ one gets
$$
\displaystyle
\limsup_{\tau\rightarrow\infty}
\frac{\displaystyle
I_{\partial\Omega}(\tau;B,\sqrt{\frac{c}{\tau}}\,)}
{\displaystyle
\tau\int_{\Omega}\vert w_0\vert^2\,dx}\le (c-1)c.
$$
This completes the proof of (1.16).

$$\quad$$

\centerline{{\bf Acknowledgment}}

The author was partially supported by Grant-in-Aid for
Scientific Research (C)(No. 17K05331) of Japan  Society for
the Promotion of Science.

$$\quad$$

\vskip1cm
\noindent
e-mail address

ikehata@hiroshima-u.ac.jp

\end{document}